\documentclass[12pt]{amsart}
\usepackage[T2A]{fontenc}
\usepackage[cp1251]{inputenc}
\usepackage{amsmath,amssymb}

\newtheorem{theorem}{Theorem}[section]
\newtheorem{lemma}{Lemma}[section]
\newtheorem{cor}{Corollary}[section]
\newtheorem{rem}{Remark}[section]
\newtheorem{defn}{Definition}[section]
\newtheorem{smt}{Proposition}[section]

\newlength\hstep
\newlength\vstep
\newcommand\col{\mathrm {col}}
\usepackage{graphicx}
\usepackage{multicol}

\title{Classification of $O$-topologically non-equivalent functions with color chord diagrams.}
\author {A.A.Kadubovsky, A.V.Klimchuk.}
\date{}
\begin{document}
\begin{abstract}
By means of color chord diagrams we establish a necessary and
sufficient condition for $O$-topological equivalence of functions
with one essentially critical point on oriented surfaces with
edge. We also calculate the number of $O$-topologically
non-equivalent functions with  one essentially critical point on
oriented surfaces with edge.
\end{abstract}

\maketitle\pagestyle{myheadings} \markboth {A.A.Kadubovsky,
A.V.Klimchuk.}{Classification of $O$-topologically equivalent
functions with color chord diagrams}
\section *{Introduction}

 Let $N$ be a finite-dimensional closed
smooth manifold, $C^\infty(N)$ be the space of infinitely
differentiable functions on $N$.

Two functions $f$ and $g\in C^\infty (N) $ with isolated critical
points are called \textit {topologically equivalent} if there
exists a homeomorphism $k:\ N \rightarrow N $ and $l:\
R^1\rightarrow R^1 $ (with $l$ preserving orientation)  such that
$g $ = $l\circ f \circ k^{-1} $.
 If $N$ is an oriented surface and the homeomorphism
  $k$ preserves orientation then functions $f$ and $g$ are called
\textit{$O$-topologically equivalent}\cite{sh2002}.

The concept of \textit {color spin-graph} is introduced
in~\cite{sh2002} and it is proved that functions $f$ and $g$ from
$C^\infty (N)$ are $O$-topologically  equivalent if and only if
there exists preserving orientation isomorphism of the color
spin-graphs.

Counting of non-isomorphic color spin-graphs with fixed Euler
characteristics is a rather difficult problem.

In this paper we introduce the notion of color chord diagram and
using these terms we give a necessary and sufficient condition of
$O$-topological equivalency for functions with one essentially
critical point on an oriented surface with edge.

We count the number of non-isomorphic $O$-diagrams and solve the
problem of calculation of $O$-topologically non-equivalent smooth
functions with one essentially critical point on oriented surfaces
with edge.

\section {Preliminaries}% Introductory definitions and notes.
Let $N$ be a smooth ($C^\infty$) compact two-dimensional manifold
with edge $\partial N$, $f$ be a smooth function on $N$, and $x\in
N$. Then the point $x\in N$ is called {\em critical} for a
function $f$ if all partial derivatives $f$ vanish at $x$.

Suppose that $x$ is an isolated critical point of $f$ being not a
local extremum of $f$.

Then there are {\em continuous} local coordinates at $x$ in which
$f=Re z^n+c$ ($n\geq 2$) provided that the topological type of
level lines in $x$ varies~\cite{Pri}. In this case $x$ will be
called {\em essential}. Otherwise $f$ reduces to the form $f=Re\
z$ and it is possible to eliminate a critical point. A critical
point of this type will be called \textit {unessential}.

Let us consider the function $f=Re z^n+c$ ($n\geq 2$, $z=x+iy$) in
a neighborhood $U$ of zero of plane $R^2$. It is obvious that the
level line $\Gamma=f^{-1}(c)$ of function $f$ in neighborhood $U$
contains the critical point  $0$ and consists of $2n$ intervals
intersecting in the point $0$ or, as we will call them below, of
$2n$ edges which are going out of one vertex. Every neighboring
pair of edges forms a sector in interior of which the function $f$
assumes value either greater or smaller than $c$ (but not equal to
$c$). Further we will call them white or black sectors
respectively. Thus $U$ includes $2n$ sequentially alternating
white and black sectors.

Let us recall some definitions and statements given in paper
\cite{sh2002}.
\begin{defn}
  Suppose that a graph $\Delta'$ consists of $2n$ edges $a_i$
  connecting vertex $x$ with $2n$ vertices $y_i$.  A color spin in
  vertex $x$, denoted $\lessdot x$, is a partition
  of edges $a_i$ into pairs $(a_i, a_j)$ together with indication of the color (black or
  white) $\col(a_i, a_j)$ for each pair, so that each edge
  makes only one white and only one black pair exactly
  with two different edges, and the color-alternating sequence $\col(a_{i_1}, a_{i_2}),
  \col(a_{i_2}, a_{i_3}), \col(a_{i_3}, a_{i_4}),\ \ldots \, \col(a_{i_{2n}}, a_{i_1})$
has the length $2n$.
\end{defn}
\begin{defn}
Suppose that the order of each vertex of graph $\Delta$ is even. A
color spin of the graph $\Delta$ is an assignment of a color spin
to each vertex of order greater than two. A graph $\Delta$, with a
color spin given,  is called a color spin-graph and is denoted as
$\lessdot\Delta$.
\end{defn}
Further if color spin-graph $\lessdot\Delta$ is given then the
color spin in the vertex $x$ is denoted as $\lessdot\Delta(x)$. It
is clear that it is possible to set a color spin on a graph in
many ways.
\begin{smt}
A graph $\Delta'$ with $2n$ edges and a color spin $\lessdot x$ in
the vertex $x$ can be imbedded in a neighborhood $U$ of the
critical point $0$ of function $Re z^n$ so that its image will lie
on $\Gamma=f^{-1} (0) \cap U$ and this imbedding will preserve
color spins into vertices at point $0$. I.e. it is possible not
only to expand graph $\Delta'$ to a disk but also to specify
function $Re z^n$ on the expansion.
\end{smt}

Let $(N, \partial_{-}N, \partial_{+}N)$ be a smooth surface with
edge $\partial N=\partial_{-}N\cup\partial_{+}N$. We will consider
levels of smooth function $f: (N,
\partial_{-}N,\partial_{+}N)\rightarrow\left[I,0,1\right]$
with only one essential critical point laying on
$f^{-1}(\frac{1}{2})$.

Then the level line $\Gamma$, containing the critical point, is
the image of imbedded into surface $N$ finite graph $\Delta$ with
one vertex of order $2n$. The vertex of the imbedded graph is the
critical point of function $f$. If we consider a neighborhood of
an essential critical point laying on $\Gamma$, in its proceeding
from function $f$ arises a color spin. Hence $\Gamma$ receives the
structure of a graph with color spin. And on the contrary, if
there is a color spin-graph it is possible to construct  a surface
according to it. Namely we expand a vertex of the color spin-graph
to a disk (proposition $1.1$). Then we color the sectors of the
disk formed by segments of edges of the spin-graph into black and
white depending on the value of a spin. Then we glue
black-and-white strips to the colored disk along graph's edges so
that the colors of sectors of the disk and the strips are
coordinated. As a result we receive a colored surface with edge
which further will be called an \textit{expansion of color
spin-graph}.

If an expansion of the color spin-graph $\lessdot\Delta$ is an
oriented surface then, having chosen its orientation, we may set
an oriented cyclical order on the edges in the neighborhood of
each vertex. In this case we speak about an \textit{oriented color
spin-graph}. If we have an isomorphism $\varphi$ between oriented
color spin-graphs which preserves bicoloured spin and orientation
in the neighborhood of each vertex we will speak that this
\textit{isomorphism $\varphi$ preserves orientation}.

Let $C_1^\infty(N,\partial N)$ be the subspace of space
$C^\infty(N)$ consisting of those functions on a surface $N$ with
edge $\partial N=\partial_{-}N\bigcup\partial_{+}N$ all of whose
critical points lie in the interior $N$ on one level line, and
they have an identical value $a$ $(b)$ on the components of
connectivity of edge $\partial_{-}N$ $(\partial_{+}N)$. In paper
\cite{sh2002} the following theorem is proved.
\begin{theorem}
There is only a finite number of $O$-topologically non-equivalent
functions from space $C_1^\infty(N, \partial N)$ on the oriented
surface $N$ with edge $\partial
N=\partial_{-}N\bigcup\partial_{+}N$ which is equal to the number
of oriented non-isomorphic color spin-graphs whose expansions are
homeomorphic to the surface $N$.
\end{theorem}

\section {Relation of color spin-graphs with color chord
diagrams.}
\begin{defn}
The configuration (actually graph) on a plane consisting of a
circle and $n$ chords connecting $2n$ various points is called a
chord diagram of the order $n$ or, shortly, an $n$-diagram
\cite{hru}.
\end{defn}
\begin{defn}
A color chord diagram is an $n$-diagram whose arcs of the circle
are colored in two colors so that any two neighboring arcs are not
colored in the same color. Denote a color chord diagram with $n$
chords by $D^\ast$ and the set of all such diagrams by
$\Im_n^\ast$.
\end{defn}
Without loss of generality it can be assumed that the points of a
color diagram have been numbered clockwise.
\begin{defn}
A color diagram which does not contain (contain) the chord
connecting points with numbers of one parity will be called an
$O$-diagram ($N$-diagram) and the set of all such diagrams is
denoted as $\Im_n^O$ ($\Im_n^N$).
\end{defn}

Let us consider the set $\Theta$ of oriented color spin-graphs
which correspond to functions with one essential critical point on
an oriented surface with edge.

Let us consider function $f:\
(N,\partial_{-}N,\partial_{+}N)\longrightarrow\left[I,0,1\right]$
and $x$ is its critical point. Let us choose an orientation on
surface $N$. Without loss of generality we may suppose that
critical points of all functions lie in one point. In the opposite
case it is not difficult to achieve this with the help of
diffeomorphism which is isotopic to identical one. A neighborhood
of critical point can be chosen so that it can be represented by a
disk neighborhood imbedded into surface with segments of lines
intersected in it (fig. $1$).
\begin{figure}[ht]
\includegraphics[width=\textwidth]{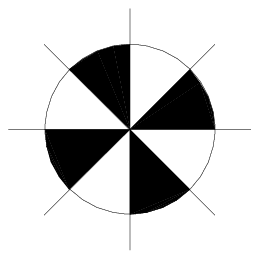}\caption{}
\end{figure}
On a critical level in a neighborhood of a critical point we set
the structure of a color spin-graph.

Then we put a color chord diagram in correspondence to the
obtained object. Namely:

On each loop $\omega_i,\ i=1,.., n$ of the color spin-graph
$\lessdot\Delta$ with a unique vertex $x$ we will put two points
$A_i, A_j$ which indexing is set according to the cyclical order
of edges already given on $\lessdot\Delta$.

Let us connect each pair of points $A_i, A_{i+1}$, $i=1,.., 2n-1$;
$A_{2n}, A_1$ by edges and color the edges according to the color
spin. Consider $$D^\ast=\left(\bigcup_{i=1}^{2n-1}A_i
A_{i+1}\bigcup A_{2n} A_1\right)
\bigcup\left(\lessdot\Delta\backslash\Lambda\right),$$ where
$\Lambda=\bigcup_{i=1}^{2n}xA_i$ is a subgraph of
$\lessdot\Delta$. Then $D^\ast$ is a color diagram whose circle is
oriented.

Thus we establish correspondence (with preservation of
orientation) of a color chord diagram $D^\ast $ to every color
spin-graph $\lessdot\Delta$ (fig. $2$).

\begin{figure}[ht]
\includegraphics[width=\textwidth]{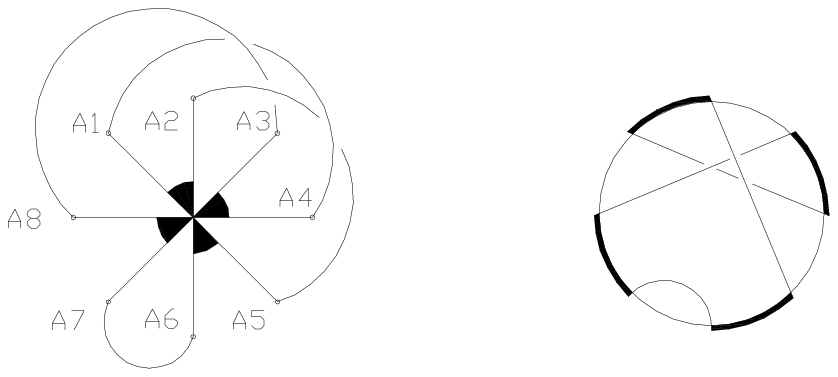}\caption{}
\end{figure}

And this correspondence $\mu: \Theta\longrightarrow\Im_n^\ast$
will be one-to-one.

Let $h:\ \lessdot\Delta_1\longrightarrow\lessdot\Delta_2$ be a
preserving orientation isomorphism of color spin-graphs. Then the
mapping $h':\ \mu\ \circ\ h\ \circ\ \mu^{-1}:\
D^\ast_1\longrightarrow D^\ast_2$ is a preserving orientation
isomorphism of oriented color diagrams. Thus if oriented color
spin-graphs are isomorphic then the corresponding oriented color
diagrams are also isomorphic. The inverse statement is correct as
well. So we obtained a bijection between the sets of color
spin-graphs and color chord diagrams and the following lemma is
proven.
\begin{lemma}
Oriented color spin-graphs are isomorphic if and only if their
corresponding color diagrams are also isomorphic.
\end{lemma}
\begin{rem}
For the functions given on oriented surfaces with edge,
corresponding oriented color spin-graphs have the following form.
On each loop going out of the point $O$ the points $A_i, A_j:
i+j=2k+1, k\in \rm{N}$ lies (see the construction of expansion of
color spin-graph). The set of $O$-diagrams corresponds to the set
of such spin-graphs.
\end{rem}

From the above the following statement is valid
\begin{lemma}
Two functions $f, g:
(N,\partial_{-}N,\partial_{+}N)\longrightarrow\left[I,0,1\right]$
with one critical point on a surface with edge are
$O$-topologically equivalent if and only if the corresponding
color chord $O$-diagrams are isomorphic.
\end{lemma}
Due to the previous lemma it is possible to reformulate the
theorem $1.1$ in terms of the color diagrams.
\begin{theorem}
There is only a finite number of $O$-topologically non-equivalent
functions with one isolated critical point on oriented surfaces.
This number is equal to number of non-isomorphic $O$-diagrams or,
what is the same, to number of the non-isomorphic color diagrams
generated by $O$-gluings \ (see def. $3.1$).
\end{theorem}
\begin{rem}
When a function is given on closed surfaces with edge and all
maxima (minima) of the function lie on one level line $L_1 (
L_2)$, the last statement is also true.
\end{rem}
\section {Gluings and color chord diagrams.}
As a pattern it is understood "colored" \ circle (the coloring is
fixed) with $2n$ sequentially numbered (clockwise) points. We will
suppose that the arcs $ \widehat{1;2},\widehat{3;4}...
\widehat{2n-1;2n}$ are black and the arcs
$\widehat{2;3},\widehat{4;5}...\widehat{2n;1}$ are white.

As a gluing $\alpha$ we will mean a way of joining $2n$ various
points of the colored  circle by $n$ chords. Any fixed point may
be contained only in one of $n$ pairs. For the greater rigor and
convenience we consider that:
\\ \centerline{$\alpha=(1=a_1, b_1)(a_2, b_2)...(a_n, b_n)\quad a_i<a_{i+1}; \
a_i<b_i,\ i=1,..,n-1$,} \\ where $a_i, b_i=1,..,2n, a_i\ne b_i$
are indices of points. A set of all such possible gluings (of $2n$
vertices of the pattern) is denoted by $\rm B_{2n}$.

It is known (\cite{hru}, \cite{Cori}) that the cardinality of the
set $\Im_n$ of all $n$-diagrams is equal to
$d_n=(2n-1)!!=\left|\rm B_{2n}\right| =\frac{(2n)!}{2^n\cdot n!}$.
Then it is obvious that the cardinality of the set of all color
chord diagrams is equal to $2\times d_n$.
\begin{smt}
The number of non-isomorphic color chord diagrams is equal to the
number of non-isomorphic diagrams (def. $4.1$) constructed on the
basis of a pattern with fixed colors.
\end{smt}
\begin{pr}
It is sufficient to show that for every diagram
$D(\alpha)\in\Im_n^\ast$ constructed on the pattern there exists a
diagram $D(\alpha')\in\Im_n^\ast$ which is isomorphic $D(\alpha)$
when the colors of $D(\alpha')$ are changed. Let
$\alpha=(a_1,b_1)...(a_n,b_n)\in\rm B_{2n}$. Then
$D(\alpha)\in\Im_n^\ast$. Consider
$\alpha'=(a_1+1\mod(2n),b_1+1\mod(2n))...(a_n+1\mod(2n),b_n+1\mod(2n))\in\rm
B_{2n}$. Then $D(\alpha')\in\Im_n^\ast$. Obvious that when we
change colors of $D(\alpha')$ and rotate on angle
$\frac{2\pi}{2n}$ (counterclockwise) then $D(\alpha)=D(\alpha')$.
\qed
\end{pr}

Therefore further we consider only color chord diagrams
constructed on the pattern. We will denote the set of all such
diagrams by $\Im_n^\ast$. Then it is natural to suppose that
\begin{equation}\label{e1}
\left|\Im_n^\ast\right|=(2n-1)!!=\frac{(2n)!}{2^n\cdot n!}
\end{equation}
\begin{defn}
The gluing $\alpha$ defining diagram $D^\ast(\alpha)$ is named
\textit{$O$-gluing} (\textit{$N$-gluing}) if it
generates an $O$($N$)-diagram.
The set of all $O$-gluings is denoted by $\rm B_{2n}^{o}$.
\end{defn}
\begin{rem}
$O$-gluings of the pattern determine oriented expansions of
corresponding color spin-graphs as far as the latter do not
contain twisted black-and-white strips.
\end{rem}
\begin{lemma}
The cardinality of the set $\Im _n^{O}$ is equal to
\begin{equation}
\label{e3} \left|\Im_n^{O}\right|=\left|\rm B_{2n}^{o}\right|=n!
\end{equation}
\end{lemma}
\begin{pr}
By definition \ $O$-gluings cannot contain the chords connecting
points of a pattern with numbers of identical parity. Thus the
$O$-gluing $\alpha$ has the form:\\ $\alpha=(1,
b_1)(3,b_2)...(2n-3, b_{n-1})(2n-1, b_n)$, where $b_i$ are even
numbers of vertices of the pattern. It is obvious that it is
possible to choose $b_1$ in $n$ ways; $b_2$ -- in $n-1$ ways;
\ldots $b_i$ -- in $i$ ways; \ldots $b_n$ -- in $1$ way. Thus
there exist exactly $n!$ \ $O$-gluings. \qed
\end{pr}
\begin{defn}
A $b$-cycle ($w$-cycle) of an expansion of color spin-graph
$\lessdot\Delta(D^\ast(\alpha))$ corresponding to the diagram
$D^\ast=D^\ast(\alpha)$ is a component of its edge colored in
black (white) color.
\end{defn}
\begin{defn}
A $b$-cycle ($w$-cycle) of diagram $D^\ast(\alpha)$ is an
alternating sequence of chords and black (white) arcs.
\end{defn}
Let us illustrate an algorithm of calculation
$\mathchar'26\mkern-10mu\lambda$ black and white cycles of diagram
on a particular example (fig. $3$).

\begin{figure}[ht]
\includegraphics[width=\textwidth]{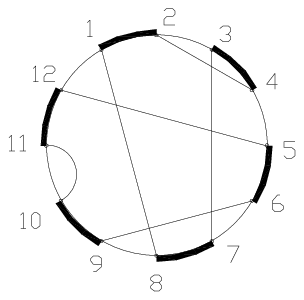}\caption{}
\end{figure}

We set an orientation on the circle of the diagram
$D^\ast=D^\ast\left(\alpha\right)$ arbitrarily (for example
clockwise). We write out all black cycles of the diagram: \\
$Cb_1=(\widehat{1,2})(2,4)(\widehat{4,3})(3,7)(\widehat{7,8})(8,1)(\widehat{1,2});$
\\ $Cb_2=(\widehat{5,6})(6,9)(\widehat{9,10})(10,11)(\widehat{11,12})(12,5)(\widehat{5,6}).$
\\ Similarly we write out all white cycles of the diagram:
\\ $Cw_1=(\widehat{2,3})(3,7)(\widehat{7,6})(6,9)(\widehat{9,8})(8,1)(\widehat{1,12})(12,5)
(\widehat{5,4})(4,2)(\widehat{2,3});$
\\ $Cw_2=(\widehat{10,11})(11,10)(\widehat{10,11})$
\\ Here $(\widehat{\cdot,\cdot})$ are arcs of the diagram; $(\cdot,\cdot)$
are its chords. It is obvious that
$\mathchar'26\mkern-10mu\lambda(b,w)=\mathchar'26\mkern-10mu\lambda(2,2)=4$.
\begin{rem}
At routing of the cycle the fixed orientation on the circle
changes direction if the diagram contains chords connecting points
with numbers of identical parity.
\end{rem}
\newpage
\section {Isomorphism of the color diagrams.}
\begin{defn}
Two diagrams $D_1^\ast(\alpha)$ and $D_2^\ast(\alpha')$ are
isomorphic, denoted $D_1^\ast(\alpha)\cong D_2^\ast(\alpha')$, if
one can be received from another as a result of some rotation.
\end{defn}
It is obvious that not colored diagram can be rotated on any angle
\\ \centerline{$\varphi=m\cdot\frac{2\pi}{2n},\  \ 0<m\le 2n$.}
\begin{rem}
It is not difficult to see that rotation of a color diagram is
correctly defined only for "even" \ angles ($m=2r$), as far as at
odd $m$ every black (white) arc of the diagram will be put over a
white (black) arc of the pattern.
\end{rem}
\begin{rem}
Under rotation of the diagram $D^\ast(\alpha)$ on angle \
$2m\cdot\frac{2\pi}{2n}\ (0<m\le n)$ we will understand such its
rotation around of the centre at which:
\\ each chord $H_i=(a_i, b_i), \ i=1,.., n$ connecting
the points with numbers $a_i, b_i$ passes in the chord
$H_i'=\left(a_i', b_i'\right)$ connecting vertices with numbers $$
a_i'=\left[
\begin{array}{ll}
[a_i+2m]\bmod(2n), & a_i+2m\ne 2n  \\ 2n, & a_i+2m=2n
\end{array}
\right.; \quad\ b_i'=\left[
\begin{array}{ll}
 [b_i+2m]\bmod(2n), & b_i+2m \ne 2n \\
 2n, & b_i+2m=2n
\end{array}
\right. : $$
\\ each arc of black color passes in an arc of black
color; the point with number\ $1$ \  passes in the point with
number $[1+2m]\bmod(2n)$; etc.
\end{rem}

Let us consider on a set $\Im_n^\ast$ of all color diagrams the
action of the subgroup \\ $G=\left\{\xi^k\in\xi: \quad k=2m, \
m=1,..,n \right\}$ of the group $\langle\xi_{2n}\rangle$ of the
order $2n$ of cyclical permutations.

In work \cite{Cori} there was established that the group
$\langle\xi_{2n}\rangle$ acts on the set $\Im_n$ as conjugation.
Namely:
\begin{smt}
Two chord diagrams $D=D(\alpha), D'=D(\alpha')$ are isomorphic if
and only if one is obtained from the other by some rotation, i.e.
if one is conjugated of the other by some power $\xi^k$ of the
circular permutation: $$D\cong
D'\Leftrightarrow\exists\xi^k\in\xi:
D(\alpha)=\xi^{-k}D(\alpha')\xi^k$$ where
$k:\quad\frac{2\pi}{2n}\cdot k$ is the angle of rotation of chord
diagram.
\end{smt}
Considering the previous proposition the following proposition is
valid.
\begin{smt}
Two color chord diagrams $D_1=D^\ast(\alpha), D_2=D^\ast(\alpha')
$ are isomorphic if and only if one is obtained from the other by
some rotation on "even" \ angle, i.e. if one is conjugated of the
other by some power $\xi^{k=2m}$ of cyclical permutation:
$$D_1\cong D_2 \Leftrightarrow\exists\ \xi^k\in G:\quad
D(\alpha)=\xi^{-k}D(\alpha')\xi^k.$$
\end{smt}
\begin{defn}
Following \cite{Cori}, a permutation $\xi^k\in G$ will be called
an \textit{automorphism} of the color diagram
$D^\ast=D^\ast(\alpha)\in\Im_n^\ast,\quad\alpha\in\rm B_{2n}$, if
$\xi^{-k}D^\ast(\alpha)\xi^k=D^\ast(\alpha)$.
\end{defn}
\begin{defn}
An \textit{orbit} $\rm O$ of a fixed diagram $D^\ast(\alpha)$
(generated by gluing $\alpha\in\rm B_{2n})$ is a subset
$\Im_n^\ast$ such that: $$\rm O=\rm
O\left(D^\ast(\alpha)\right)=\{D^\ast(\beta)\in\Im_n^\ast\ |\
D^\ast(\beta)=\xi^{-k}D^\ast(\alpha)\xi^k, \quad \xi^k\in G \}. $$
\end{defn}
\begin{defn}
A \textit{stabilizer} $G_{D^\ast(\alpha)}$ of a chord diagram
$D^\ast(\alpha)$ is a subset of $G$ such that:
$$G_D^\ast(\alpha)=\left\{\xi^k\in G:\
D^\ast(\alpha)=\xi^{-k}D^\ast(\alpha)\xi^k\right\}.$$
\end{defn}
\begin{defn}
By
$\chi\left(\xi^k,n,\Upsilon\right)=\chi\left(\xi^k,n\right)_\Upsilon$
let us denote the set of fixed diagrams \\
$D^\ast(\beta)\in\Upsilon\subseteq\Im_n^\ast$ of the cyclical
permutation $\xi^k\in G: $
$$\chi(\xi^k,n)_\Upsilon=\left\{D^\ast(\beta)\in\Upsilon:
D^\ast(\beta)=\xi^{-k}D^\ast(\beta)\xi^k,\ \xi^k\in G\right\}.$$
\end{defn}
On the set $\rm B_{2n}$ of all gluings let us define the operation
$\Re(\alpha,k)$ of "rotation on angle" \\ $1\le k\le 2n $ \ as
follows: $\Re(\alpha,k)=\alpha+k\bmod(2n)=\alpha'$.
\\ I.e.: $\alpha+k\bmod(2n)=$
\\ $=\left([1+k]\bmod(2n),[b_1+k]\bmod(2n)\right)\ ...\
\left([a_n+k]\bmod(2n),[b_n+k]\bmod(2n)\right)$, where: $$
d+k\bmod(2n)=\left[\begin{array}{ll}
 d+k\bmod (2n), & d+k\bmod (2n) \ne 0 \\
 2n, & d + k\bmod (2n)=0
  \end{array}
  \right.,\quad
d=\left[\begin{array}{ll}
 a_i \\
 b_i
 \end{array}
 \right.
$$ Then it is not difficult to see that:
\begin{equation}
\label{e5}
\Re(\alpha,k)=\alpha+k\bmod(2n)=\alpha'\Leftrightarrow\xi^{-k}D^\ast(\alpha)\xi^k=D^\ast(\alpha').
\end{equation}
\begin{theorem}
For any $n\ge 2$ the number $d_n^{\ast\ast}$ of non-isomorphic
color chord diagrams is calculated by the formula:
\begin{equation}
\label{e6} d_n^{\ast\ast}=\frac{1}{n}\left[(2n-1)!!\ +
\sum\limits_{k\vert2n,\ k=2m\ne
2n}{\phi\left(\frac{n}{m}\right)\rho(n,2m)}\right],
\end{equation}
where: $(2n - 1)!!=\left|\Im_n^{\ast}\right|$; $\phi(q)$ is Euler
arithmetic function (the number of integers smaller than $q$ which
are relatively prime with it);
\begin{equation} \label{e7}
p(n,k)=p(n,2m)=\chi\left(\xi^k,n\right)_{\Upsilon=\Im
_n^\ast}=\left[\begin{array}{ll}
 (2m - 1)!!\cdot\left(\frac{n}{m}\right)^m, & \frac{n}{m}=2l+1 \\
 \sum\limits_{r = 0}^m {C_{2m}^{2r}\cdot (2r - 1)!!
 \cdot\left(\frac{n}{m} \right)^r}, & \frac{n}{m}=2l
 \end{array} \right.
 \end{equation}
\end{theorem}
\begin{pr}
From lemma $3.1$ \cite{Cori} it follows that the number of
non-isomorphic diagrams may be calculated by the formula:
\begin{equation}\label{e8}
d_n^\star=\frac{1}{2n}\left((2n-1)!!+\sum\limits_{\xi^k\ne\
\xi^{2n}\in G: \ k\vert
2n}{\phi\left(\frac{2n}{k}\right)Fix\left(\xi^k,n\right)}\right),
\ \mbox{where:}
\end{equation}
\\ $2n=\left|\langle\xi_{2n}\rangle\right|$;\quad
$Fix\left(\xi^k,n\right)=\left|\left\{D(\alpha )\in\Im _n:
\xi^{-k}D(\alpha)\xi^k=D(\alpha)\right\}\right|=$\\
$=\left|\left\{ \alpha\in{\rm B}_{2n}: \Re(\alpha, k)=\alpha +
k\bmod (2n)=\alpha\right\}\right|$ is the number of diagrams for
which the permutation $\xi^k\in\langle\xi_{2n}\rangle$ is an
automorphism.
\begin{equation}
\label{e9}
 Fix\left(\xi^k,n\right)=p(n,k)=
 \left[\begin{array}{ll}
 (k - 1)!!\cdot\left(\frac{2n}{k}\right)^{\frac{k}{2}}, & \frac{2n}{k}=2l+1 \\
 \sum\limits_{r=0}^{\left[k/2\right]}{C_k^{2r}\cdot(2r - 1)!!
\cdot\left(\frac{2n}{k}\right)^r }, & \frac{2n}{k}=2l
\end{array}
\right.
\end{equation}
By virtue \textit{proposition $4.2$}, from the relations
$(\ref{e8})$ and $(\ref{e9})$ follows that
$$d_n^{\ast\ast}=\frac{1}{n}\left((2n - 1)!!\ +
\sum\limits_{k=2m\ne 2n, \ k\vert 2n} {\phi
\left(\frac{n}{m}\right)\cdot p(n,k)}\right), \ \mbox{where}$$
$p(n,k)=p(n,2m)=\left[ {\begin{array}{ll}
 (2m - 1)!! \cdot \left( {\frac{n}{m}} \right)^m, & \frac{n}{m}= 2l+1 \\
 \sum\limits_{r = 0}^m {C_{2m}^{2r} \cdot (2r - 1)!!\cdot\left(\frac{n}{m}\right)^r}, & \frac{n}{m} = 2l
 \end{array}} \right.$\qed
\end{pr}
\begin{cor}
For any prime $n\ge3$ the number of non-isomorphic color chord
diagrams is equal to:
\begin{equation}
\label{e10}
d_n^{\ast\ast}=\frac{(2n-1)!!}{n}+n-1=\frac{(2n)!}{2^n\cdot
n!\cdot n}+n-1
\end{equation}
\end{cor}
\begin{pr}
From \textit{theorem $4.1$} it follows that the number of
non-isomorphic color diagrams is equal to
 $d_n^{\ast\ast}=\frac{1}{n}\left((2n - 1)!!\ +\sum\limits_{\xi^{k=2m}\ne\ \xi^{2n}\ \in G}
 \phi\left(\frac{2n}{2m}\right)\cdot p(n, 2m)\right).$
As among even divisors of the number $2n$ there are only $2$ and
$2n$ then $d_n^{\ast\ast}=\frac{1}{n}\left[(2n-1)!!\ +
\phi\left(n\right)\cdot p(n, 2)\right]$. It is obvious that
$\phi\left(n\right)=n-1,\ p(n, 2)=p(n, 2\cdot 1)=(2 \cdot 1-1)!!\
\times\left(\frac{n}{1}\right)^1=n.$
\\ Therefore $d_n^{\ast\ast}=\frac{1}{n}\left[(2n-1)!!+(n-1)n\right].$\qed
\end{pr}

Calculation of non-isomorphic color chord diagrams for $n=2,3...,
11$ gives the following values.
\begin{table}[ht]
\begin{tabular}
{|p{88pt}|p{88pt}|p{88pt}|p{88pt}|} \hline $n$& $d_n^{\ast \ast }
$& $n$& $d_n^{\ast \ast } $ \\ \hline 2& 3& 7& 19311 \\ \hline 3&
7& 8& 254143 \\ \hline 4& 35& 9& 3828921 \\ \hline 5& 193& 10&
65486307 \\ \hline 6& 1799& 11& 1249937335 \\ \hline
\end{tabular}
\label{tab1}
\end{table}
\section {Non-isomorphic color diagrams generated by $O$-gluings}
It was established earlier that
\begin{itemize}
  \item any $O$-gluing has the form
$\alpha=(1,b_1)(3,b_2)...(2i-1,b_i)...(2n-1,b_n)$, \\ where
$b_i\in\left\{2,4..., 2n\right\} \ \forall i=1,.., n; \ b_i\ne
b_j$;
  \item the number of the color diagrams generated by all possible
$O$-gluings, \\ is equal to $n!$, \ i.e., the cardinality of set
$\Im_n^O$; \item  the subgroup
$G=\left\{\xi^k\in\langle\xi_{2n}\rangle: \ k=2m, \
m=1,..,n\right\}$ of group of cyclical permutations acts on the
set $\Im_n^O$; \item the permutation $\xi^k\in G$ is an
automorphism of the diagram $D^\ast=D^\ast(\alpha)\in \Im_n^O,$ \\
$\alpha\in\rm B_{2n}^O$, if $\xi^{-k}D(\alpha)\xi^k=D(\alpha)$;
\item
$\Re(\alpha,k)=\alpha+k\bmod(2n)=\alpha'\Leftrightarrow\xi^{-k}D(\alpha)\xi^k
=D'(\alpha'). $
\end{itemize}
Then from lemma  $3.1$ \cite{Cori} and theorem $4.1$ it follows
that the number of non-isomorphic  color diagrams generated by
$O$-gluings is calculated by the formula:
\begin{equation} \label{e11}
d_n^{\ast(O)}=\frac{1}{n}\left(n! +
\sum\limits_{\xi^k\ne\xi^{2n}\in G:\
k=2m\vert2n}{\phi\left(\frac{2n}{k}\right)\chi(\xi^k,n)_{\Upsilon=\Im_n^O}}\right),\
\hbox{where:}
\end{equation}
$n!=\left|\Im_n^O\right|$ is the cardinality of the set of color
diagrams generated by $O$-gluings;
\\ $\phi(q)$ is Euler arithmetic function;
 \quad $\chi\left(\xi^k, n\right)_{\Im_n^O}=\left|\left\{D(\alpha)\in\Im_n^{O}:\ \xi^{-k}
 D(\alpha)\xi^k=D(\alpha)\right\}\right|=$
 \\ $=\left|\left\{\alpha\in\rm B_{2n}^O: \Re(\alpha,k)=\alpha+k\bmod(2n)=\alpha\right\}\right|$
 is the number of diagrams for which the permutation $\xi^k$ is an automorphism.

\begin{lemma}
For any $n=i\cdot m$
\begin{equation}\label{e12}
\chi\left(\xi ^{2i},n\right)_{\Im_n^O}=\prod\limits_{l=1}^i{l\cdot
m}=(i)! \cdot m^i=(i)! \cdot \left(\frac{n}{i}\right)^i
\end{equation}
In particular $\forall n\geq 2\in{\rm N}$
\begin{equation}\label{e13}
   \chi\left(\xi^2, n\right)_{\Im_n^O}=n
   \end{equation}
\end{lemma}
\begin{pr}
The proof is reduced to the calculation of the number of all
\textit{$O$-gluings} \ $(i)! \cdot m^i\alpha$ satisfying the
condition: $\Re(\alpha, k=2m)=\alpha+k\bmod(2n)=\alpha$. All such
gluings have the form \\
$\alpha=(1,b_1)(3,b_2)(5,b_3)...(2i-1,b_i)...(2n-1,b_n)$, where
$b_i\in\left\{2,4,.., 2n\right\};\ \forall i=1..., n; b_i\ne b_j$.
\\ $k=2:\ $ As far as $\Re(\alpha,2)=\alpha+k\bmod(2n)=\alpha$ then
$b_i=b_{i-1}+2$. But then all these gluings have the form: \\
$\alpha=(1,b_1)(3,b_1+2)(5,b_1+4)...(2i-1,b_1+2i-2)...(2n-1,b_1+2n-2)$.
It is obvious that the number of such gluings is equal to $n$. The
latter proves $(\ref{e13})$.
\\ $k=4:\ $ As far as $\Re(\alpha,4)=\alpha+k\bmod(2n)=\alpha$, then
$b_i=\ b_{i-2}+4$. But then all such gluings decompose into two
subgluings $\alpha_1, \alpha_2\ : \alpha=\alpha_1;\alpha_2$ which
have the form:
\\ $\alpha=(1,b_1')(5,b_1'+4)(9,b_1'+8)(13,b_1'+12)...(4i-3,b_1'+4i-4)...
(2n-3,b_1'+2n-4)$;
\\ $(3,b_2')(7,b_2'+4)(11,b_2'+8)(15,b_2'+12)...(4j-1,b_2'+4j-4)...(2n-1,b_2'+2n-4)$.

As far as $4|2n$ then $n=2m$. It is not difficult to see that in
the subgluings $i, j=1,..,m$. But then the number of all such
gluings is equal to $2m\cdot(n-m)=2m\cdot m$.
\\ $k=6:\ $ As far as $6|2n$ then $n=3m$. As far as $\Re(\alpha,6)=\alpha$,
then $b_i=b_{i-3}+6$. But then all such gluings decompose into
three subgluings $\alpha_1, \alpha_2, \alpha_3:
\alpha=\alpha_1;\alpha_2;\alpha_3 $  which also have the form:
\\ $\alpha=(1,b_1')(7,b_1'+6)(13,b_1'+12)(19,b_1'+18)...(6i-5,b_1'+6i-6)...(2n-5,b_1'+2n-6)$;
\\ $(3,b_2')(9,b_2'+6)(15,b_2'+12)(21,b_2'+18)...(6j-3,b_2'+6j-6)...(2n-3,b_2'+2n-6)$;
\\ $(5,b_3')(11,b_3'+6)(17,b_3'+12)(23,b_3'+18)...(6r - 1,b_3'+6r-6)...(2n-1,b_3'+2n-6)$.

It is not difficult to see that in subgluings $i,j,r=1,..,m$. But
then the number of all such gluings is equal to
$3m\cdot(n-m)((n-m)-m)=3m\cdot 2m\cdot m$.
\\ \ldots
\\ $k=2i:\ $ As far as $2i|2n$ then
$n=i\cdot m $. As far as $\Re(\alpha,2i)=\alpha$ then all gluings
decompose into $i=\frac{n}{m}$ subgluings with $m$ elements in
each. Then repeating the reasoning similar to one used in first
two cases, we obtain that the number of all such gluings is equal
to $im\cdot(im-m)\cdot((im-m)-m)\cdot...\cdot
m=\prod\limits_{l=1}^i{l\cdot m}$ what proves the validity
$(\ref{e12})$.
\end{pr} \qed
\begin{cor}
The number of non-isomorphic color $O$-diagrams is calculated by
formula:
\begin{equation}\label{nf1}
d_n^{\ast(O)}=\frac{1}{n}\left(n!\ + \sum\limits_{i\vert n\ i\ne
n}{\phi\left(\frac{n}{i}\right)\cdot
i!\cdot\left(\frac{n}{i}\right)^i}\right)
\end{equation}
\end{cor}
\begin{cor}
Let $p\geq3$ be a prime number. Then
\begin{equation} \label{14}
  d_p^{\ast(O)}=(p-1)!+(p-1)
  \end{equation}
\end{cor}
\begin{cor}
The number of non-isomorphic color $N$-diagrams is calculated by
formula:
\begin{equation}
\label{nf2} d_n^{\ast(N)}=d_n^{\ast\ast)} - d_n^{\ast(O)}
\end{equation}
\end{cor}
\begin{table}[htbp]
\begin{tabular}
{|p{88pt}|p{88pt}|p{88pt}|p{88pt}|} \hline $n$& $d_n^{\ast(O)}$ &
$n$& $d_n^{\ast(O)} $ \\ \hline 2& 2& 7& 726 \\ \hline 3& 4& 8&
5100 \\ \hline 4& 10& 9& 40362 \\ \hline 5& 28& 10& 363288
\\ \hline 6& 136& 11& 3628810 \\ \hline
\end{tabular}
\end{table}
{\bfseries Acknowledgement.} The authors are grateful to Vladimir
Sharko for posing the problem and to Andrey Khruzin for
stimulating conversations.


\begin{thebibliography}{99}
\bibitem{sh2002} Шарко\ В.В. Гладкая и топологическая эквивалентность функций
на поверхностях. // Укр. мат. жур. - 2003. - 55, №\ 5 - С.687-700.
\bibitem{Bol} Болсинов~А.В., Фоменко~Ф.Т. Введение в топологию
 интегрируемых гамильтоновых систем. М.: Наука 1997-352 с.
\bibitem{Pri} Prishlyak A.O.
Topological equivalence of smooth functions with isolated critical
points on a cloused surface.
//Topology and its Aplications.-2002.-119.-p.257-267.
\bibitem{Man} Мантуров~В.О. Атомы, высотные атомы, хордовые диаграммы и
 узлы. Перечисление атомов малой сложности с использованием языка
 Mathematica 3.0 // Топологические методы в теории гамильтоновых
 систем (сборник статей) под редакцией А.В. Болсинова, А.Т.Фоменко,
 А.И. Шафаревича. М. Изд-во  Факториал - 1998, С. 203-212.
\bibitem{Sto} Stoimenov~A. On the number of chord diagrams.// Discrete
Math.- 2000.-218. N1-3-p.209-233.
\bibitem{Cori} R. Cori, M. Marcus.  Counting non-isomorphic chord
diagrams. // Teoretical Computer Science - 1998-204. -p. 55-73.
\bibitem{hru} A. Khruzin. Enumeration of chord diagrams. - Arxiv: math. CO/0008209,10p.
\end{thebibliography}
\end{document}